\documentclass{article}
\usepackage[utf8]{inputenc}
\usepackage{graphicx}
\usepackage{psfrag}
\usepackage{epsf}
\usepackage{amsmath,amsfonts,amssymb,latexsym,amsthm, comment,amsthm, enumerate}
\usepackage{algorithmic}
\usepackage{algorithm}
\usepackage[margin=1.1in]{geometry}
\usepackage{setspace}
\usepackage{enumerate}
\usepackage{tikz}
\usepackage{soul}
\usepackage{xcolor}
\usepackage{hyperref}

\usetikzlibrary{calc}
\tikzstyle{line}=[draw]
\usetikzlibrary{shapes.geometric}

\newcommand{\ignore}[1]{}
\newcommand{\startClaims}{\setcounter{claim}{0}}
\newtheorem{theorem}{Theorem}[section]
\newtheorem{corollary}[theorem]{Corollary}
\newtheorem{porism}[theorem]{Porism}
\newtheorem{lemma}[theorem]{Lemma}
\newtheorem{proposition}[theorem]{Proposition}
\newtheorem{conjecture}[theorem]{Conjecture}

\newtheorem{definition}[theorem]{Definition}

\newtheorem{observation}[theorem]{Observation}

\newtheorem*{corollary*}{Corollary}

\theoremstyle{plain}

\newcommand{\eopf}{\raisebox{0.8ex}{\framebox{}}}

{\noindent{\bf Sketch of a proof.}\ }%
{\hfill\eopf\par\bigskip}%
{\noindent{\bf Ideas for a proof.}\ }%
{\hfill\eopf\par\bigskip}%
{\noindent{\bf Proof in progress.}\ }%
{\hfill\eopf\par\bigskip}%

\newcommand{\jdup}{\operatorname{jdup}}

\newcommand{\Sym}{\mathcal{S}}

\newcommand{\Gbar}{\overline{G}}

\title{Graphs with Bipartite Complement that Admit Two Distinct Eigenvalues}

\author{Wayne Barrett
\thanks{Department of Mathematics, Brigham Young University, Provo, UT (wb@mathematics.byu.edu) }
\and Shaun Fallat
\thanks{Department of Mathematics and Statistics, University of Regina, Regina, Saskatchewan, CA (shaun.fallat@uregina.ca)} 
\and Veronika Furst \thanks{Department of Mathematics, Fort Lewis College, Durango, CO (furst\_v@fortlewis.edu)}
\and Shahla Nasserasr \thanks{School of Mathematical Sciences, Rochester Institute of Technology, Rochester, NY (sxnsma@rit.edu)} ~\thanks{Corresponding Author}
\and Brendan Rooney \thanks{School of Mathematical Sciences, Rochester Institute of Technology, Rochester, NY (brsma@rit.edu)} 
\and Michael Tait \thanks{Department of Mathematics \& Statistics, Villanova University, Villanova, PA (michael.tait@villanova.edu) }}

\date{\today}

\numberwithin{equation}{subsection}

\begin{document}
\maketitle

\begin{abstract}
The parameter $q(G)$ of an $n$-vertex graph $G$ is the minimum number of distinct eigenvalues over the family of symmetric matrices described by $G$. We show that all $G$ with $e(\overline{G}) \leq \lfloor n/2 \rfloor -1$ have $q(G)=2$. We conjecture that any $G$ with $e(\overline{G}) \leq n-3$ satisfies $q(G) = 2$.   We show that this conjecture is true if $\overline{G}$ is bipartite and in other sporadic cases. Furthermore, we characterize $G$ with $\overline{G}$ bipartite and $e(\overline{G}) = n-2$ for which $q(G) > 2$. 
\end{abstract}

\noindent {\bf Keywords} inverse eigenvalue problem for graphs, orthogonality, $q$-parameter, strong spectral property. \\

\noindent {\bf AMS subject classification} 05C50, 15A29, 15A18. \\

\section{Introduction} \label{intro}

Let $G = (V(G), E(G))$ be a simple graph on $n = |V(G)|$ vertices and $e(G) = |E(G)|$ edges.  Let $\Sym(G)$ be the set of $n\times n$ real symmetric matrices for which $a_{ij} = 0$, $i\neq j$, if and only if $ij \notin E(G)$.  The Inverse Eigenvalue Problem for Graphs (IEPG) asks which spectra can occur for matrices in $\Sym(G)$ \cite{MR4478249}.  One approach to the IEPG is to investigate the minimum number of distinct eigenvalues that can occur for a matrix $A\in\Sym(G)$ \cite{MR3118943}.  Given an $n\times n$ symmetric matrix $A$, let $q(A)$ denote the number of distinct eigenvalues of $A$ and let 
\[q(G) = \min\{q(A)\,:\,A\in\Sym(G)\}. \]

The determination of $q(G)$ is a difficult outstanding problem for almost all graphs, but it has been achieved for $n\leq 6$ \cite{MR3904092}; see also \cite{MR4284782}.  A compelling subproblem is determining which graphs have $q(G) = 2$ (since only the graph with no edges satisfies $q(G) = 1$); see for example \cite{MR4044603, MR3118943, MR3904092, MR3506498, MR3891770, MR4774297}. One approach is to find a lower bound on the number of edges that a graph on $n$ vertices with $q(G) = 2$ can have.  The solution to this so-called ``allows problem" was given in \cite{MR4603836}:  If a connected graph $G$ on $n$ vertices satisfies $q(G) = 2$, then $e(G)\geq 2n-4$ if $n$ is even and $e(G)\geq 2n-3$ if $n$ is odd.  Equality can be achieved for any $n$ and graphs achieving equality were characterized in \cite{MR4603836}.  Another approach to the $q(G) = 2$ problem is to narrow the focus to a subclass of graphs, such as regular graphs.  The classification of regular graphs of degree at most 4 with $q(G)=2$ appears in \cite{MR4646344}.

In this paper, we shift our focus from graphs of low density to graphs of high density. The $q=2$ ``requires problem" poses the following query:  For a fixed number of vertices, $n$, determine the largest positive integer $m$ such that the removal of any set of $m$ edges from the complete graph $K_n$ results in a graph $G$ that admits a matrix $A\in\Sym(G)$ with $q(A)=2$.  For $n\geq 2$, it is well-known that complete graphs all have $q(K_n)=2$ (Lemma 2.2 in \cite{MR3118943}). Furthermore, if $n\geq 4$ and $e$ is any edge of $K_n$ then $q(K_n-e)=2$ (Proposition 4.7 in \cite{MR3118943}).  For $n=5$ and $n=6$, the minimum number of edges needed for a graph to require $q=2$ is $8$ and $12$, respectively \cite{MR3904092}. Equivalently, from the standpoint of the complement, a graph $G$ on 5 vertices with $e(\Gbar)\leq 2$ admits a matrix $A\in\Sym(G)$ with $q(A)=2$, as does a graph $G$ on 6 vertices with $e(\Gbar)\leq 3$.  It therefore seems reasonable to pose the following conjecture (also listed explicitly as Conjecture 1 in \cite{FallatMojallal}), which has been verified for $n=7$ and $n=8$ \cite{FallatMojallal}.

\begin{conjecture} \label{mainconj}
    If $G$ is a graph on $n\geq 3$ vertices such that $e(\Gbar) \leq n-3$, then $q(G) = 2$.
\end{conjecture}

We require a few additional definitions and background results, the first of which imposes a very useful necessary, though not sufficient, condition for $q(G)=2$.

\begin{lemma}[Corollary of \cite{MR3118943}, Theorem 3.2] \label{nouniquenb}
 Let $G$ be a connected graph with $q(G) = 2$.  If $xuy$ is a path of length 2, then either $x$ and $y$ are adjacent or there exists vertex $v\neq u$ such that $xvy$ is another path of length 2.
\end{lemma}

We note that Conjecture \ref{mainconj} cannot be improved further, as demonstrated by the following observation.

\begin{observation}
  Removing $n-2$ edges incident to a single vertex in $K_n$ ($n\geq 3$) results in a connected graph with a vertex of degree 1.  By Lemma \ref{nouniquenb}, this graph cannot have $q$-value of 2.
\end{observation}

Given a vertex $v\in V(G)$, let $N(v)$ denote the set of neighbors (adjacent vertices) of $v$ and let $N[v] = N(v) \cup \{v\}$.  Both notations can be extended to a set $S\subseteq V(G)$, and we may use a subscript to indicate the relevant graph when necessary.  The open neighborhood of $v$ has cardinality $\deg(v) = |N(v)|$, and a subscript may again be used to denote the relevant graph.  A graph $\jdup(G,v)$ is constructed from $G$ by \emph{joined-duplicating} a vertex $v\in V(G)$ if $V(\jdup(G,v))=V(G)\cup\{u\}$ and $E(\jdup(G,v))=E(G)\cup\{uw\,:\,w\in N[v]\}$.

\begin{lemma}[\cite{MR3891770}, Lemma 2.9] \label{jdup}
For any vertex $v$ in a connected graph $G$, $q(\jdup(G,v))\leq q(G)$.
\end{lemma}

Let $G\vee H$ denote the {\em join} of $G$ and $H$, the graph consisting of the vertices $V(G)\cup V(H)$ and all edges in $E(G) \cup E(H) \cup \{ uv\,:\,u\in V(G), \ v\in V(H)\}$.

\begin{proposition}[\cite{MR3506498}, Theorem 5.2] \label{joinSameOrder}
    For any two connected graphs $G$ and $H$ with the same number of vertices, $q(G\vee H) = 2$.
\end{proposition}

\begin{proposition}[\cite{MR4774297}, Theorem 3.4] \label{joinExtension}
    For any two connected graphs $G$ and $H$ whose orders differ by at most 2, $q(G\vee H) = 2$.
\end{proposition}

The following lemma is a generalization of the two previous propositions in a different direction.  The technique in the proof highlights the Strong Spectral Property \cite{MR3665573}, which we apply frequently throughout this paper.  A symmetric matrix $A$ is said to have the {\em Strong Spectral Property (SSP)} if the only symmetric matrix $X$ (known as the certificate matrix) satisfying $A\circ X = O$, $I\circ X = O$, and $AX - XA = O$ is $X = O$.  Here $\circ$ denotes the Hadamard (entrywise) product and $O$ denotes the zero matrix.  Theorem 10 of \cite{MR3665573} states that if $A\in\Sym(G)$ has the SSP, then every supergraph $G'$ of $G$ with the same vertex set has a realization $A'\in\Sym(G')$ such that $A'$ has the same spectrum as $A$ and has the SSP.

\begin{lemma}\label{JoinClique}
If $H$ is any graph on $n$ vertices and $n'\geq n$, then $q(K_{n'}\vee H) = 2$.
\end{lemma}

\begin{proof} 
Let $G=K_n$ in Lemma 3.13 in~\cite{MR4044603}. Then $q(K_n\vee \overline{K_n})=2$. In addition, the $(1,2)$ block of the matrix $M\in \mathcal{S}(K_n\vee \overline{K_n})$ in the proof of the lemma is invertible. Now, Lemma 5.1 of \cite{JOS19} implies that $M$ has the SSP. The result follows from joined-duplication of the vertices of $K_n$. 
\end{proof}
 
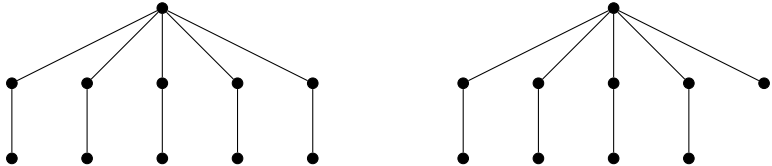
\begin{figure}[ht]
\begin{center}
\begin{tikzpicture}[scale=1, vrtx/.style args = {#1/#2}{%
      circle, draw, fill=black, inner sep=0pt,
      minimum size=4pt, label=#1:#2}]
\node (b1) [vrtx=right/] at (6,0) {};
\node (b2) [vrtx=right/] at (4,-1) {};
\node (b3) [vrtx=above/] at (5,-1) {};
\node (b4) [vrtx=above/] at (6,-1) {};
\node (b5) [vrtx=left/] at (7,-1) {};
\node (b6) [vrtx=right/] at (8,-1)  {};
\node (b7) [vrtx=right/] at (4,-2) {};
\node (b8) [vrtx=right/] at (5,-2) {};
\node (b9) [vrtx=right/] at (6,-2) {};
\node (b10) [vrtx=right/] at (7,-2) {};
\draw (b1) edge (b2);
\draw (b1) edge (b3);
\draw (b1) edge (b4);
\draw (b1) edge (b5);
\draw (b1) edge (b6);
\draw (b2) edge (b7);
\draw (b3) edge (b8);
\draw (b4) edge (b9);
\draw (b5) edge (b10);

\node (a1) [vrtx=right/] at (0,0) {};
\node (a2) [vrtx=right/] at (-2,-1) {};
\node (a3) [vrtx=above/] at (-1,-1) {};
\node (a4) [vrtx=above/] at (0,-1) {};
\node (a5) [vrtx=left/] at (1,-1) {};
\node (a6) [vrtx=right/] at (2,-1)  {};
\node (a7) [vrtx=right/] at (-2,-2) {};
\node (a8) [vrtx=right/] at (-1,-2) {};
\node (a9) [vrtx=right/] at (0,-2) {};
\node (a10) [vrtx=right/] at (1,-2) {};
\node (a11) [vrtx=right/] at (2,-2) {};
\draw (a1) edge (a2);
\draw (a1) edge (a3);
\draw (a1) edge (a4);
\draw (a1) edge (a5);
\draw (a1) edge (a6);
\draw (a2) edge (a7);
\draw (a3) edge (a8);
\draw (a4) edge (a9);
\draw (a5) edge (a10);
\draw (a6) edge (a11);
\end{tikzpicture}
\caption{The graphs $W(5,0,\vec{1})$ and $W(4,1,\vec{1})$.}\label{fig:Wgraphs}
\end{center}
\end{figure}

Throughout this paper we denote the subgraph of $G$ induced by $S\subseteq V(G)$ as $G[S]$.  The notation $H \subseteq G$ indicates that $H$ is a subgraph of $G$ (equivalently, $G$ is a supergraph of $H$) and $V(H) = V(G)$.  The removal of an edge $e$ (as above) or a vertex $v$ from a graph $G$ are denoted $G-e$ and $G-v$, respectively.  The Cartesian product of two graphs $G$ and $H$ is written as $G\square H$.  Using the notation of \cite{MR3891770}, for $a,b\geq 0$ let $S_{a,b}$ denote a tree whose vertex set can be partitioned as $\{x_1,\ldots,x_a\}\cup \{u,v\}\cup \{y_1,\ldots,y_b\}$ so that each $x_i$ is a leaf adjacent to $u$, each $y_i$ is a leaf adjacent to $v$, and $u$ is adjacent to $v$.  When $a,b\geq 1$, $S_{a,b}$ is known as a \emph{double-star}.  As in \cite{MR3891770}, for $\vec{1} = (1,1, \ldots, 1)$, let $W(k,0,\vec{1})$ and $W(k,1,\vec{1})$ be the graphs obtained from $K_{1,k}$ by subdividing all edges or from $K_{1,k+1}$ by subdividing all edges except one, respectively (see Figure \ref{fig:Wgraphs}).

This paper is organized as follows.  In Section \ref{Requires} we present some results towards the full generality of Conjecture \ref{mainconj} but with fewer removed edges.  In Section \ref{BipartiteComplement} we establish our principal result (Theorem \ref{n-3 edges}) that the conjecture holds under the additional assumption that $\Gbar$ is bipartite. Moreover, for $\Gbar$ bipartite and $e(\Gbar) = n-2$, we establish in Section \ref{Characterization} (Theorem \ref{tightbipartitechar}) that $q(G) = 2$ unless $\Gbar$ has a specific structure, in which case $q(G) = 3$. We have an ever increasing quantity of evidence that the bipartite restriction can be removed, some of which is included in Section \ref{OtherComplements}.

\section{The \texorpdfstring{$q=2$}{q=2} Requires Problem}\label{Requires}

In this section, we prove a bound on $e(\overline{G})$ that ensures $q(G)=2$ (Theorem \ref{LB1}). For general graphs this is the closest result we have towards Conjecture \ref{mainconj}.

\begin{observation}
    If $G$ is a graph on $n\geq 3$ vertices such that $e(\Gbar) \leq \lfloor n/4 \rfloor$, then $q(G) = 2$.
\end{observation}

\begin{proof}
   The statement holds for $n\leq 7$ by Lemma 2.2 and Proposition 4.7 in \cite{MR3118943}.  For $n>7$, the removal of at most $\lfloor n/4 \rfloor $ edges from $K_n$ results in 
   $K_{\lceil n/2 \rceil} \vee H$ where $H$ is a connected graph on $\lfloor n/2 \rfloor$ vertices, and we may apply Proposition \ref{joinExtension}.
\end{proof}

With some effort, using Lemma \ref{JoinClique}, the upper bound in the previous observation can be improved to $e(\Gbar) \leq \lfloor 3n/8 \rfloor$.  Our best general answer to the ``requires problem" is given in Theorem \ref{LB1}, whose proof makes use of the following combinatorial lemma.

\begin{lemma} \label{lemmaT}
    Suppose $t_1\geq t_2\geq \cdots\geq t_k$ are natural numbers with $\sum_{i=1}^k t_i = n$ and $\sum_{i=1}^k (t_i - 1) \leq \frac{n}{2} - 1$.  Then there exist disjoint sets $A,B$ with $A\cup B = \{1, \ldots, k\}$ such that \[
    \left| \sum_{i\in A} t_i - \sum_{j\in B} t_j \right| \leq 1.
    \]
\end{lemma}

\begin{proof}
By the assumptions we note that $k\geq \frac{n}{2}+1$. Let $x$ be the number of $i$ such that $t_i = 1$ and $y$ be the number of $i$ such that $t_i>2$. Our first claim is that $x\geq y+2$. This follows because 
\[
\frac{n}{2}-1 \geq \sum_{i=1}^k (t_i - 1) \geq 2y + 1(k-y-x) + 0x = y+k-x \geq y+\frac{n}{2}+1 - x;
\]
rearranging gives the claim. 

We will prove the statement by induction on $t_1$. If $t_1 = 1$ or $t_1=2$, the claim gives that $t_k = t_{k-1} = 1$. Let $m$ be the number of $t_i$ that are $2$; that is, $t_1=\cdots = t_m = 2$ and $t_{m+1} = \cdots = t_k = 1$. If $m$ is even or $k$ is even, let $A = \{t_{2i+1}\,:\,0\leq i \leq \lfloor (k-1)/2\rfloor \}$ and $B = \{t_{2i}\,:\,1\leq i \leq \lfloor k/2\rfloor \}$. If $m$ is odd and $k$ is odd, let $A = \{t_1, t_3,\ldots t_{k-4}, t_{k-2}\}$ and $B = \{t_2, t_4, \ldots, t_{k-3}, t_{k-1}, t_k\}$.
In either case, $A$ and $B$ give a satisfactory partition.

Now assume that the statement holds for sequences with greatest value strictly less than $l\geq 3$ and satisfying the hypotheses of the statement. Suppose we have $t_1,\ldots, t_k$ satisfying the hypotheses of the lemma and $t_1=l$. If $k > x+y$, consider the sequence of $k-y$ integers given by 
\[
t_i' = \begin{cases}
    t_i - 1 & \mbox{ \ for \ } 1\leq i \leq y\\
    2 & \mbox{ \ for \ } y+1 \leq i \leq k-x\\
    1 & \mbox{ \ for \ } k-x+1 \leq i \leq k-y;
\end{cases}
\]
otherwise, let 
\[
t_i' = \begin{cases}
    t_i - 1 & \mbox{ \ for \ } 1\leq i \leq y\\
    1 & \mbox{ \ for \ } k-x+1 \leq i \leq k-y.
\end{cases}
\]

Note that the last range is well-defined by the claim. This sequence $t_1',\ldots, t_{k-y}'$ is non-increasing, satisfies $t_1' = t_1-1<l$, and
\[
\sum_{i=1}^{k-y} (t_i' - 1) = \left(\sum_{i=1}^k (t_i - 1)\right) - y \leq \frac{n}{2}-1 - y = \frac{n-2y}{2}-1 = \frac{\sum_{i=1}^{k-y}t_i'}{2}-1,
\]
and we may therefore apply the induction hypothesis. So there exists a partition $A'$ and $B'$ such that 
\[
\left| \sum_{i\in A'} t_i' - \sum_{j\in B'} t_j'\right| \leq 1.
\]

Start with sets $A$ and $B$ where $t_i \in A$ if and only if $t_i' \in A'$ and $t_j \in B$ if and only if $t_j' \in B'$. Then 
\[
\left| \sum_{i\in A} t_i - \sum_{j\in B} t_j\right| \leq 1 + y,
\]
and the integers $t_{k-y+1},\cdots ,t_k$ are all equal to $1$ and have not been assigned to either $A$ or $B$ yet. Since they are all equal to $1$ it is easy to assign each to $A$ or $B$ such that
\[
\left| \sum_{i\in A} t_i - \sum_{j\in B} t_j\right| \leq 1.
\]

\end{proof}

\begin{theorem}\label{LB1}
    If $G$ is a graph on $n\geq 3$ vertices such that $e(\Gbar) \leq \lfloor n/2 \rfloor - 1$, then $q(G) = 2$.
\end{theorem}

\begin{proof}
Assume $G$ is a graph with $e(G) \geq \binom{n}{2} - \lfloor \frac{n}{2} \rfloor + 1$ edges. We will first assume $n$ is even.  Let $H_1, \ldots, H_k$ be the connected components of $\overline{G}$ with $|V(H_i)| = t_i$.  Without loss of generality, assume $t_1 \geq t_2 \geq \cdots \geq t_k$.  

Note that 
\begin{equation} \label{t-1ineq}
\sum_{i=1}^k (t_i - 1) \leq \sum_{i=1}^k e(H_i) \leq \left\lfloor \frac{n}{2} \right\rfloor - 1.
\end{equation}
By Lemma \ref{lemmaT}, there exist sets $A,B$ such that $A\cap B = \emptyset$, $A\cup B = \{1, \ldots, k\}$, and
\[
    \left| \sum_{i\in A} t_i - \sum_{j\in B} t_j \right| \leq 1.
\]

Since $n$ is even we have that $\sum_{i\in A} t_i = \sum_{j\in B} t_j$.
By the first claim in the proof of Lemma \ref{lemmaT}, the components $H_{k-1}$ and $H_k$ are isolated vertices, so they are dominating vertices in $G$.  If $\{k-1, k\} \not\subseteq A$ and $\{k-1, k\} \not\subseteq B$, then \[
 G_A=G\left[\bigcup_{i\in A} V(H_i)\right] \quad\text{and}\quad G_B=G\left[\bigcup_{j\in B} V(H_j)\right]
\]
are connected subgraphs of $G$. Since $G=G_A\vee G_B$, and $G_A$ and $G_B$ have the same order, by Proposition \ref{joinSameOrder}, we have that $q(G) = 2$.  If one of the sets in the partition, say $B$, contains both $k-1$ and $k$, then we may move $k$ to $A$, so that $G_A$ and $G_B$ are both connected; then $G = G_A \vee G_B$ with $|G_A| - |G_B| = 2$, and $q(G) = 2$ follows by Proposition \ref{joinExtension}.

Finally, assume that $n$ is odd. Again, since $H_{k-1}$ and $H_k$ are isolated vertices, they are dominating vertices in $G$.  In fact, $H_{k-2}$ must also be an isolated vertex.  Suppose otherwise.  Then $n = \sum t_i \geq 2(k-2) + 2$, which implies $n/2 \geq k-1$.  On the other hand, it follows from Inequality \ref{t-1ineq} that $\lceil \frac{n}{2} \rceil \leq k-1$.  Since $n/2 < \lceil n/2 \rceil$, this gives a contradiction.  Let $H_k = \{v\}$, $H_{k-1} = \{w\}$, and define $G' = G - v$.  Since $e(G) \geq \binom{n}{2} - \frac{n}{2} + \frac{3}{2}$, we see that 
\[ 
e(G') = e(G) - (n-1) \geq \binom{n}{2} - \frac{n}{2} + \frac{3}{2} - (n-1) =  \binom{n-1}{2} - \left\lfloor \frac{n-1}{2}\right\rfloor + 1.
\]
By the case when $n$ is even, we have that $q(G') =2$. Since $G = \jdup(G', w)$, Lemma \ref{jdup} now implies $q(G) = 2$.     
\end{proof}

\section{Graphs With Bipartite Complement}\label{BipartiteComplement}
 In this section, we prove that Conjecture \ref{mainconj} is true when $\overline{G}$ is bipartite (Theorem \ref{n-3 edges}). We first give some definitions and lemmas.   For a large subclass of graphs whose complement is bipartite, we use the following theorem of Levene et al. \cite{MR3891770} to establish Conjecture \ref{mainconj}.  Recall that $q(K_n) = 2$, so we only consider graphs with nonempty complement.

\begin{theorem}[\cite{MR3891770}, Theorem 2.7]\label{NGThm}
Let $\overline{G}$ be a bipartite graph with partite sets of sizes $r,s \geq 1$, where $G$ has order $r+s \geq 3$. If $\overline{G}$ does not contain $K_{m_1,n_1} \cup K_{m_2,n_2}$ as a spanning subgraph, where $r = m_1+m_2$, $s = n_1+n_2$, $m_1, m_2, n_1, n_2 \geq 0$, and either $m_1m_2 \neq 0$ or $n_1n_2 \neq 0$, then $q(G) = 2$. 
\end{theorem}

\begin{definition}\label{SimplifiedDefinition}
    A graph $G$ on $m_1+m_2+n_1+n_2\geq 3$ vertices is called {\em simplified} if $\overline{G}$ is nonempty and bipartite, with partite sets $M_1\cup M_2\neq \emptyset$ and $N_1\cup N_2\neq\emptyset$, such that $M_1\cup N_1$ induces $K_{m_1, n_1}$ and $M_2\cup N_2$ induces $K_{m_2, n_2}$ in $\overline{G}$, and furthermore there is no pair of non-isolated vertices in $\overline{G}$ with exactly the same neighborhood.
\end{definition}

To motivate the previous definition, we make the following observation, which we will use repeatedly in our arguments, allowing us to restrict our attention to simplified graphs.
\begin{observation}\label{whysimplified}
Let $G$ be a graph on at least three vertices (with a nonempty complement).
    If $N_{\overline{G}}(u) = N_{\overline{G}}(v)\neq\emptyset$ for vertices $u$ and $v$, then $\jdup(G-v, u) = G$; note that $e(\overline{G-v})\leq e(\Gbar) - 1$. We sequentially remove one vertex from each pair of non-isolated vertices with exactly the same neighborhood in $\overline{G}$ until we obtain a graph $G'$ in which no such pairs exist.  If $e(\overline{G})\leq k$, then $e(\overline{G'}) \leq k$.  Furthermore, $q(G) \leq q(G')$ by Lemma \ref{jdup}, and in particular if $q(G')=2$, then $q(G)=2$ .
\end{observation}

\begin{lemma}\label{simplified properties}
    Let $G$ be a simplified graph of order $n\geq 3$ such that $e(\overline{G}) \leq n-2$ or $e(\overline{G}) = n-1$ and $\overline{G}$ contains a cycle. Partition the vertices of $\Gbar$ as in Definition \ref{SimplifiedDefinition}. If $G$ is not $C_4$, then the following hold:
    \begin{enumerate}
        \item At least one of $m_1, m_2, n_1, n_2$ equals 0. Without loss of generality, assume that $N_2$ is empty.
        \item It follows from $n_2 = 0$ that $m_1 = 1$.
        \item If $n$ is even, then $n_1 \leq \frac{n}{2}-1$.
        \item If $n$ is even, then there are at least $\frac{n}{2}-n_1$ isolated vertices in $\overline{G}$.
    \end{enumerate}
\end{lemma}

\begin{proof}
To establish the first assertion, we consider the possibilities where none of $m_1$, $m_2$, $n_1$, or $n_2$ is $0$. First suppose that $m_1, n_1, m_2, n_2\geq 2$. Then
\[
n-1\geq e(\overline{G})\geq m_1n_1+m_2n_2\geq m_1+n_1+m_2+n_2=n,
\]
giving an immediate contradiction. Without loss of generality, assume $m_1=1$. 
First suppose that $m_2=1$. Then $\overline{G}[M_1\cup N_1]$ and $\overline{G}[M_2\cup N_2]$ are both stars and account for $n-2$ edges. If $e(\overline{G})=n-2$ then $n_1\geq 2$ or $n_2\geq 2$ contradicts our assumption that $G$ is simplified, and $n_1=n_2=1$ contradicts $G\neq C_4$. If $e(\overline{G})=n-1$, then our assumption that $\overline{G}$ contains a cycle implies that $\overline{G}$ contains a triangle, contradicting that $\Gbar$ is bipartite. This same argument gives a contradiction when $n_2=1$.

Finally, suppose $m_2, n_2\geq 2$. We have
\[
n-1\geq e(\overline{G})\geq m_1n_1+m_2n_2\geq n_1+m_2+n_2=n-1,
\]
so we have equality throughout.  This implies that $m_2=n_2=2$ and $\overline{G}= K_{1,n_1}\cup C_4$. However this means the vertices in $M_2$ (and the vertices in $N_2$) have the same neighborhood, contradicting $G$ being simplified.

Since we have ruled out every possibility where none of $m_1$, $m_2$, $n_1$, or $n_2$ is $0$, one of them must be. We assume without loss of generality that $n_2=0$.

Since $N_2 = \emptyset$, all $\overline{G}$-neighbors of vertices in $M_1$ belong to $N_1$.  Since $M_1\cup N_1$ induces $K_{m_1,n_1}$ in $\overline{G}$, it follows from $G$ being a simplified graph that $m_1 = 1$, establishing assertion (2).

For assertion (3), note that $N_1$ can contain at most one vertex with $\overline{G}$-degree 1 since $G$ is simplified.  So 
\[ e(\overline{G}) \geq n_1 + n_1 - 1.\]
If $e(\overline{G})\leq n-2$, then the claim follows since $n$ is even and $n_1$ is an integer. Suppose $e(\overline{G}) = n-1$ and $\overline{G}$ contains a cycle; then $n_1 \leq n/2$ follows immediately.  If $n_1 = n/2$, then $N_1$ must contain one vertex with $\overline{G}$-degree 1 and all remaining vertices with $\overline{G}$-degree 2; but this implies that the cycle in $\overline{G}$ has order $4$ and that the two cycle vertices in $N_1$ have the same neighborhood, a contradiction.

Finally, let $I$ represent the set of isolated vertices in $\overline{G}$ (note that $I\subseteq M_2$), and let 
\[ E^1 = \{e \in E(\overline{G})\,:\,e \mbox{ is incident to a vertex } v\in M_2 \mbox{ with } \deg_{\overline{G}}(v) = 1\}, \]
\[ E^2 = \{e \in E(\overline{G})\,:\,e \mbox{ is incident to a vertex } v\in M_2 \mbox{ with } \deg_{\overline{G}}(v) > 1\}, \]
\[ V^1 = \{ v\in M_2\,:\,\deg_{\overline{G}}(v) = 1\},\]
and
\[ V^2 = \{ v\in M_2\,:\,\deg_{\overline{G}}(v) > 1\}.\]
Then
\[n_1 + m_2  = n - 1 \geq e(\overline{G}) = n_1 + |E^1| + |E^2|.\]
On the other hand, $m_2 = |I| + |V^1| + |V^2|$, so 
\[ |I| + |V^1| + |V^2|  \geq |E^1| + |E^2| = |V^1| + |E^2|.\]
Since $|E^2| \geq 2|V^2|$, the previous inequality implies
\[ |I| \geq |E^2| - |V^2|  \geq |V^2|  = m_2 - |I| - |V^1| .\]

Note that $|V^1| \leq n_1$:  indeed, if $|V^1| > n_1$, then a pair of vertices in $M_2$, each with $\overline{G}$-degree 1, share their $\overline{G}$-neighbor in $N_1$, violating $G$ being simplified.  This fact together with the previous inequality yields
\[ 2|I| \geq m_2 - |V^1|  = n - n_1 - |V^1| -1 \geq n - 2n_1-1,\]
and assertion (4) follows since $n$ is even and $|I|$ is an integer.
\end{proof}

Note that $q(C_4) = 2$ is known, so the previous lemma is not weakened by excluding this graph from consideration.

\begin{lemma}
\label{boxproduct}
    Let $n\geq 4$ be even, and let $G$ be a simplified graph of order $n$ such that either $e(\overline{G}) \leq n-2$ or $e(\overline{G}) = n-1$ and $\overline{G}$ contains a cycle. Then $G$ is a supergraph of $K_{n/2} \square K_2$.
\end{lemma}

\begin{proof}
It is easy to verify that when $n=4$, $G = K_4-e$ or $G= C_4$, and the conclusion holds.  Suppose $n\geq 6$.  By Lemma \ref{simplified properties}, we may assume that we have arranged the isolated vertices of $\overline{G}$ so that $\overline{G}$ is a bipartite graph with partite sets $X$ and $Y$ that each contain $n/2$ vertices.  We will use Hall's Theorem to show the existence of a perfect matching in $G$ between $X$ and $Y$.  Let $H$ be the bipartite subgraph of $G$ with vertex set $V(H)=V(G) = X\cup Y$ and edge set $E(H)=\{uv\in E(G)\,:\,u\in X, v\in Y\}$.  Suppose $S$ is a subset of $X$ such that $|S| = k$ and $|N_H(S)| < k$. 
Let $s\in S$.  Then $\deg_H(s) \leq k-1$, so $\deg_{\overline{G}}(s) \geq \frac{n}{2}-k + 1$.  It follows that 
\[ n-1 \geq e(\overline{G}) \geq k\left( \frac{n}{2} - k + 1 \right). \]
    Rearranging, we have 
    \[
    (k+1)(k-2) + 1 = k^2-k-1 \geq \frac{kn}{2} - n = \frac{n}{2}(k-2).
    \]
If $k>2$, we divide by $(k-2)$ and obtain $n/2 \leq k+1 + \frac{1}{k-2}$.  If $k>3$, this implies $k\geq \frac{n}{2}-1$ since $k$ is an integer. Therefore, we may assume that either $k\leq 3$ or $k\geq \frac{n}{2}-1$. 

Suppose $k=3$.  Then $n\leq 10$.  The three vertices in $S$ each must be adjacent in $\overline{G}$ to all the vertices in $Y\setminus N_H(S)$, and $|Y\setminus N_H(S)|\geq n/2-2$. If $n=10$, then $e(\overline{G})\leq 9$ implies $|N_H(S)|=2$ and $\overline{G}[S\cup (Y\setminus N_H(S))]= K_{3,3}$. This accounts for all of the edges of $\overline{G}$ and contradicts that $G$ is simplified. If $n=8$, then $e(\overline{G})\leq 7$ implies $|N_H(S)|=2$ and $\overline{G}[S\cup (Y\setminus N_H(S))]= K_{3,2}$. This accounts for six of the edges in $\overline{G}$, and the addition of any edge cannot eliminate the existence of two non-isolated vertices with the same neighborhood in $\overline{G}$, again contradicting that $G$ is simplified. The final case to consider is $n=6$.  Again $|N_H(S)| = 2$, and up to isomorphism, only the configurations of edges in $\overline{G}$ shown in Figure \ref{configurations} are candidates for $e(\overline{G}) \in \{3, 4, 5\}$.

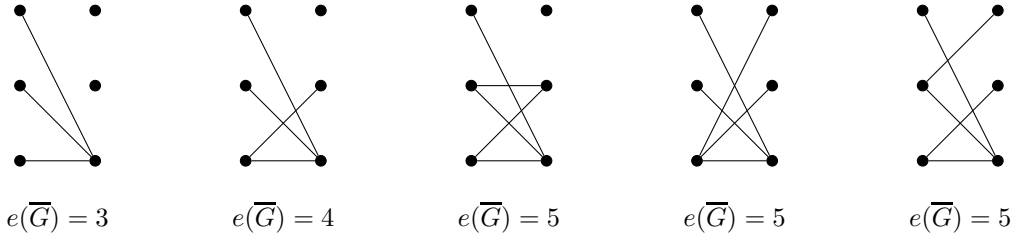
\begin{figure}[ht]
\centering
\begin{tikzpicture}[scale=1, vrtx/.style args = {#1/#2}{%
      circle, draw, fill=black, inner sep=0pt,
      minimum size=4pt, label=#1:#2}]

\node (a1) [vrtx=left/] at (-3,-1) {};
\node (a2) [vrtx=right/] at (-3,0) {};
\node (a3) [vrtx=above/] at (-3,1) {};
\node (a4) [vrtx=above/] at (-2,-1) {};
\node (a5) [vrtx=left/] at (-2,0) {};
\node (a6) [vrtx=right/] at (-2,1)  {};
\node (l1) [] at (-2.5,-1.75) {$e(\overline{G})=3$};
\draw (a1) edge (a4);
\draw (a2) edge (a4);
\draw (a3) edge (a4);


\node (b1) [vrtx=left/] at (0,-1) {};
\node (b2) [vrtx=right/] at (0,0) {};
\node (b3) [vrtx=above/] at (0,1) {};
\node (b4) [vrtx=above/] at (1,-1) {};
\node (b5) [vrtx=left/] at (1,0) {};
\node (b6) [vrtx=right/] at (1,1)  {};
\node (l2) [] at (0.5,-1.75) {$e(\overline{G})=4$};
\draw (b1) edge (b4);
\draw (b1) edge (b5);
\draw (b2) edge (b4);
\draw (b3) edge (b4);


\node (d1) [vrtx=left/] at (3,-1) {};
\node (d2) [vrtx=right/] at (3,0) {};
\node (d3) [vrtx=above/] at (3,1) {};
\node (d4) [vrtx=above/] at (4,-1) {};
\node (d5) [vrtx=left/] at (4,0) {};
\node (d6) [vrtx=right/] at (4,1)  {};
\node (l4) [] at (3.5,-1.75) {$e(\overline{G})=5$};
\draw (d1) edge (d4);
\draw (d1) edge (d5);
\draw (d2) edge (d4);
\draw (d2) edge (d5);
\draw (d3) edge (d4);


\node (e1) [vrtx=left/] at (6,-1) {};
\node (e2) [vrtx=right/] at (6,0) {};
\node (e3) [vrtx=above/] at (6,1) {};
\node (e4) [vrtx=above/] at (7,-1) {};
\node (e5) [vrtx=left/] at (7,0) {};
\node (e6) [vrtx=right/] at (7,1)  {};
\node (l5) [] at (6.5,-1.75) {$e(\overline{G})=5$};
\draw (e1) edge (e4);
\draw (e1) edge (e5);
\draw (e1) edge (e6);
\draw (e2) edge (e4);
\draw (e3) edge (e4);


\node (c1) [vrtx=left/] at (9,-1) {};
\node (c2) [vrtx=right/] at (9,0) {};
\node (c3) [vrtx=above/] at (9,1) {};
\node (c4) [vrtx=above/] at (10,-1) {};
\node (c5) [vrtx=left/] at (10,0) {};
\node (c6) [vrtx=right/] at (10,1)  {};
\node (l3) [] at (9.5,-1.75) {$e(\overline{G})=5$};
\draw (c1) edge (c4);
\draw (c1) edge (c5);
\draw (c2) edge (c4);
\draw (c2) edge (c6);
\draw (c3) edge (c4);


\end{tikzpicture}
\caption{The five candidate configurations for $e(\overline{G})\in\{3,4,5\}$.}\label{configurations}
\end{figure}
\noindent In the first four cases, $G$ is not simplified while in the last two cases, $\overline{G}$ has no cycle. So $k=3$ is impossible.

If $k=1$ or $k = n/2$, then $|N_H(S)| < k$ implies there exists a vertex $v$ (in $X$ or $Y$, respectively) with $\deg_{\overline{G}}(v) = n/2$.  If $k=2$ or $k = \frac{n}{2} - 1$, then $|N_H(S)| < k$ implies that two vertices (in $X$ or $Y$, respectively) share $\frac{n}{2}-1$ $\overline{G}$-neighbors; since $G$ is simplified, one of these vertices, $v$, must have $\deg_{\overline{G}}(v) = n/2$.  Suppose, without loss of generality, that $v\in X$.  Since $G$ is simplified, for every $y\in Y$ except at most one, there exists $x\in X$ such that $x\neq v$ and $xy \in E(\overline{G})$.  This accounts for at least $\frac{n}{2}-1$ additional edges, and 
\[  n-1 \geq e(\overline{G}) \geq \frac{n}{2} + \frac{n}{2} - 1.\]
Thus we have equality throughout, and $e(\overline{G}) = n-1$. So $Y$ must consist of one vertex with $\overline{G}$-degree 1 and $\frac{n}{2}-1$ vertices with $\overline{G}$-degree 2, each pair of which share only $v$ as their common neighbor in $X$ (since $G$ is simplified).  Thus $\overline{G}$ is the tree on $n$ vertices obtained from the star $K_{1,n/2}$ by subdividing all but one edge.  Note that this graph is $W(\frac{n}{2}-1, 1, \vec{1})$. As $\overline{G}$ has no cycle, we reach a contradiction, and $|N_H(S)| \geq |S|$.

Since $|X| = |Y|$, there exists a perfect matching $M$ in $H$ by Hall's Theorem.  It follows that the graph with vertex set $V(G)$ and edge set $E(G[X]) \cup E(G[Y]) \cup M$ is $K_{n/2} \square K_2$.  Therefore, $G$ is a supergraph of $K_{n/2} \square K_2$.
 
\end{proof}

\begin{lemma} \label{isolated}
Let $G$ be a simplified graph of order $n \geq 4$ such that $e(\overline{G}) \leq n-3$ or $e(\overline{G}) = n-2$ and $\overline{G}$ contains a cycle. Then $\overline{G}$ has at least $2$ isolated vertices.
\end{lemma}

\begin{proof}
Because $G$ cannot be $C_4$, we may apply Lemma \ref{simplified properties}. Using the notation established in the statement and proof of Lemma \ref{simplified properties}, we have $m_1=1$ and $n_2=0$. Let $M_2$ be split into vertices which have degree $0$, $1$, or at least $2$ in $\overline{G}$ and, as before, call these sets $I, V^1, V^2$, respectively. Then 
\begin{align*}
e(\overline{G}) &\geq n_1 + 2|V^2| + |V^1| \ \mbox{ and } \\
n & = n_1 + 1 + |I| + |V^1| + |V^2|.
\end{align*}
If $e(\overline{G}) \leq n-3$, it follows that $|I| \geq |V^2| + 2 \geq 2$.  On the other hand, if $e(\overline{G}) = n-2$ and $\overline{G}$ contains a cycle, then $|V^2| \geq 1$ and $|I| \geq |V^2| + 1 \geq 2$.
\end{proof}

\begin{theorem} \label{n-3 edges}
    Let $G$ be a graph of order $n\geq 3$ such that $\overline{G}$ is bipartite and $e(\overline{G}) \leq n-3$.  Then $q(G) = 2$.
\end{theorem}

\begin{proof}
As noted in Section \ref{intro}, $q(G) = 2$ if $e(\Gbar) = 0$ or $e(\Gbar) = 1$, so we may assume $n\geq 5$ (and $\Gbar$ is nonempty).  If the hypotheses of Theorem \ref{NGThm} are met, then the result follows, so assume that there is a bipartition of $\overline{G}$ with partite sets $M_1\cup M_2$ and $N_1\cup N_2$, both nonempty, such that $M_1\cup N_1$ induces $K_{m_1, n_1}$ and $M_2\cup N_2$ induces $K_{m_2, n_2}$ in $\overline{G}$.  By Observation \ref{whysimplified}, we may assume that $G$ is simplified.

If $n$  is even, Theorem \ref{boxproduct} implies $G$ is a supergraph of $K_{n/2} \square K_2$.  Suppose $n$ is odd.  By Lemma \ref{isolated}, the number of isolated vertices in $\overline{G}$ is at least 2.  Let $w$ and $z$ be isolated vertices in $\overline{G}$, and let $G' = G - w$.  Then $G = \jdup(G',z)$, $|G'| = n' = n-1$ is even, and $e(\overline{G'}) \leq n'-2$.  It suffices to show that $q(G') = 2$ by Lemma \ref{jdup}.  By Theorem \ref{boxproduct}, $G'$ is a supergraph of $K_{n'/2} \square K_2$.  By Theorem 20 of \cite{FallatMojallal}, $K_s \square K_2$, for any $s\geq 3$, has an SSP matrix realization with 2 distinct eigenvalues; the same is true when $s=2$ \cite{MR4074182}.  Therefore, any supergraph of this graph will have $q$-value 2. 
\end{proof}

\section{Graphs Close To The Requires Bound}\label{Characterization}

A natural next step is to consider those graphs that are missing one more edge above the bound in Theorem \ref{n-3 edges}, distinguishing graphs that have a $q$-value of 2 from graphs that do not. Recall that $S_{a,b}$ is the (possibly degenerate) double star with $a+b$ leaves. In this section, we prove the following characterization of graphs $G$ with $\overline{G}$ bipartite and $e(\overline{G})=n-2$ in terms of their $q$-values. 

\begin{theorem}\label{tightbipartitechar} Let $G$ be a graph of order $n\geq 3$.
If $\overline{G}$ is bipartite and $e(\overline{G})=n-2$, then $q(G)=3$ if and only if $\overline{G}= S_{a,b}\cup K_1$ for some $a,b\geq 0$, and $q(G)=2$ otherwise.
\end{theorem}

The remainder of this section provides a proof of Theorem \ref{tightbipartitechar}. First note that if $G$ meets the hypotheses of Theorem \ref{NGThm}, $q(G) = 2$ follows immediately.  So we suppose otherwise.  By Observation \ref{whysimplified}, we may assume that $G$ is simplified and $e(\Gbar) = n-2$.  (Note that if the removal of non-isolated vertices with the same neighborhood yields a simplified graph $G$ that satisfies $e(\Gbar) < n-2$, then $q(G) = 2$ follows by Theorem \ref{n-3 edges}.)

Recall that if $n\geq 4$ is even, Lemma \ref{boxproduct} implies that $G$ is a supergraph of $K_{n/2}\square K_2$. In this case, as in the proof of Theorem \ref{n-3 edges}, Theorem 20 of \cite{FallatMojallal} implies that $q(G)=2$.

Next, suppose $n$ is odd and $\overline{G}$ contains a cycle. From Lemma \ref{isolated} if $\overline{G}$ contains a cycle, then it must have at least $2$ isolated vertices. Let $w$ and $z$ be isolated in $\overline{G}$. Consider $G'=G-w$ and note that $G=\jdup(G',z)$, $|G'|=n'=n-1$ is even, and $e(\overline{G'})=n'-1$. So by Lemma \ref{boxproduct}, $G'$ is a supergraph of $K_{n'/2}\square K_2$, and Theorem 20 of \cite{FallatMojallal} implies $q(G')=2$. Therefore $q(G)=2$.

Finally, suppose that $n$ is odd and $\overline{G}$ does not contain a cycle. If $n=3$, then $G= P_3$ and $q(G)=3$ and Theorem \ref{tightbipartitechar} holds ($\overline{P_3}=S_{0,0}\cup K_1$).

For $n\geq 5$ we have the following characterization of simplified graphs $G$ that satisfy the assumptions of Theorem \ref{tightbipartitechar}.
\begin{lemma} \label{W-tree union K_1}
    Let $G$ be a simplified graph of odd order $n\geq 5$, such that $e(\overline{G}) = n-2$ and $\overline{G}$ does not contain a cycle. Then $\overline{G} = W(k, 1, \vec{1}) \cup K_1$ where $k = (n-3)/2$.
\end{lemma}
\begin{proof}
    First we apply Lemma \ref{simplified properties}. Since $n$ is odd, $G\neq C_4$. So we can partition the vertices as in Definition \ref{SimplifiedDefinition}, and conclude that $n_2=0$ and $m_1=1$.

    Since $\overline{G}$ does not contain a cycle, every vertex in $M_2$ has degree at most $1$ in $\overline{G}$. Partition $M_2$ as $M_2=I\cup V^1$ where $I$ consists of the vertices with degree $0$ in $\overline{G}$ and $V^1$ consists of the vertices with degree $1$ in $\overline{G}$. Since $G$ is simplified, any pair of vertices in $V^1$ has distinct neighbors in $\overline{G}$. Thus each vertex in $V^1$ is matched to a distinct vertex in $N_1$.
    
    Since $G$ is simplified, there is at most one vertex in $N_1$ with degree $1$ in $\overline{G}$. So we have that $|V^1|\in\{n_1-1,n_1\}$. Now from $n = 1+n_1+|V^1|+|I|$ and $n-2 = e(\overline{G})=n_1+|V^1|$, it follows that $|I|=1$. Since $n\geq 5$ is odd, we must have $|V^1|=n_1-1>0$ and so $\overline{G}= W(|V^1|, 1, \vec{1}) \cup K_1$.
\end{proof}

If $k=1$ in Lemma \ref{W-tree union K_1}, then $\overline{G}= P_4\cup K_1$ and $q(G)=3$ (see graph G47 in \cite{MR4284782}), so Theorem \ref{tightbipartitechar} holds ($P_4\cup K_1=S_{1,1}\cup K_1$).

From Lemma \ref{W-tree union K_1}, to complete the proof of Theorem \ref{tightbipartitechar}, under the condition that $G$ is simplified, it suffices to show $q(G)=2$ for all graphs $G$ with $\overline{G}=W(k, 1, \vec{1}) \cup K_1$ with $k\geq 2$. Our proof will use the Strong Spectral Property for all such values of $k$. For simplicity we denote the complement of $W(k,1,\vec{1})\cup K_1$ by $G_{2k+3}$. We consider the graph $G_7$ (i.e., $k=2$) first and verify that  $q(G_{7})=2$.

\begin{figure}[ht]
\begin{center}
\begin{tikzpicture}[scale=1, vrtx/.style args = {#1/#2}{%
      circle, draw, fill=black, inner sep=0pt,
      minimum size=4pt, label=#1:#2}]

\node (a1) [vrtx=right/] at (0,1) {};
\node (a2) [vrtx=right/] at (0,-1) {};
\node (a3) [vrtx=above/] at (1,0) {};
\node (a4) [vrtx=above/] at (2,-1) {};
\node (a5) [vrtx=left/] at (2,0) {};
\node (a6) [vrtx=right/] at (2,1)  {};
\node (a7) [vrtx=right/] at (3,0) {};
\draw (a1) edge (a2);
\draw (a1) edge (a3);
\draw (a2) edge (a3);
\draw (a1) edge (a6);
\draw (a2) edge (a4);
\draw (a3) edge (a5);
\draw (a4) edge (a7);
\draw (a5) edge (a7);
\draw (a6) edge (a7);
\end{tikzpicture}
\caption{The graph $\overline{H_7}$.}\label{H7}
\end{center}
\end{figure}
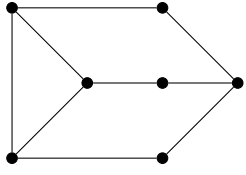

\begin{lemma}\label{W-graphs,k=2}
$q(G_{7})=2$.
\end{lemma}
\begin{proof}
Consider the matrix
\[
M_7 = \left[\begin{array}{rrrrrrr}
3 & 1 & 1 & 1 & 0 & 0 & 0 \\
1 & 1 & 0 & 0 & 1 & -1 & 0 \\
1 & 0 & 1 & 0 & -1 & 0 & 1 \\
1 & 0 & 0 & 1 & 0 & 1 & -1 \\
0 & 1 & -1 & 0 & 2 & -1 & -1 \\
0 & -1 & 0 & 1 & -1 & 2 & -1 \\
0 & 0 & 1 & -1 & -1 & -1 & 2
\end{array}\right].
\]
It is easy to check that $q(M_7)=2$ and that $M_7$ has the SSP. The off-diagonal zero pattern of $M_7$ corresponds to the graph $H_7$ whose complement is shown in Figure \ref{H7}. Since $G_7$ is a supergraph of $H_7$, hence it follows that $q(G_7)=2$.

%
%
%
%
\end{proof}

For $k\geq 3$ we use the same strategy as in Lemma \ref{W-graphs,k=2}, but we will use a different subgraph of $G_{2k+3}$. Denote by $W_{2k+3}$ the graph $W(k+1,0,\vec{1})$. Since $W(k,1,\vec{1})\cup K_1\subseteq W_{2k+3}$, we have $\overline{W_{2k+3}}\subseteq G_{2k+3}$. Our goal is to construct a matrix $M_{2k+3}\in\mathcal{S}(\overline{W_{2k+3}})$ with $q(M_{2k+3}) = 2$ and the SSP. We begin with a construction and observations from \cite{MR3240026} (see Theorem 2.4).

\begin{lemma}\label{construct}
For any $k\geq 3$, there is a matrix $B$ and vector $\vec{v}$ that satisfy the following:
\begin{enumerate}
\item $B\in\mathcal{S}(K_{k+1})$,
\item the spectral radius of $B$ is strictly less than $1$,
\item there is a totally nonzero unit vector $\vec{v}$ in the null space of $B$,
\item $B^2$ is entrywise positive.
\end{enumerate}
\end{lemma}
\begin{proof}
Let $T_{k+1}$ be the $(k+1)\times (k+1)$ matrix with $[T_{k+1}]_{i,j}=(i-j)^2$ for all $1\leq i,j\leq k+1$. It follows easily that $T_{k+1}\in\mathcal{S}(K_{k+1})$ and $T_{k+1}$ has zero diagonal. Since the entries of $T_{k+1}$ are all non-negative and for $i\neq j$ they are positive, it follows that $T_{k+1}^2$ is entrywise positive.

We use the observation from the proof of Theorem 2.4 in \cite{MR3240026} that if $\vec{r}_i$ is the $i$th row of $T_{k+1}$, then for all $4\leq i\leq k+1$,
\[
\vec{r}_i=\frac{1}{2}(i-2)(i-3)\vec{r}_1-(i-1)(i-3)\vec{r}_2+\frac{1}{2}(i-1)(i-2)\vec{r}_3.
\]
Now for all $4\leq i\leq k+1$, let $\vec{u}_i$ be the vector defined as
\[
[\vec{u}_i]_j=\begin{cases}
\frac{1}{2}(i-2)(i-3) & \mbox{ for } j=1\\
-(i-1)(i-3) & \mbox{ for } j=2\\
\frac{1}{2}(i-1)(i-2) & \mbox{ for } j=3\\
-1 & \mbox{ for } j=i\\
0 & \text{ else}.
\end{cases}
\]
Then $T_{k+1}\vec{u}_i=\vec{0}$ for each $4\leq i\leq k+1$. Define $\vec{u}=\sum_{i=4}^{k+1} \vec{u}_i$. Then $T_{k+1}\vec{u}=\vec{0}$, and $\vec{u}$ is a totally nonzero vector (to see this note that $[\vec{u}]_j=-1$ for all $4\leq j\leq k+1$, and when $j\in\{1,2,3\}$, $[\vec{u}]_j$ is either the sum of $k-2$ positive values or $k-2$ negative values).

Now since $T_{k+1}$ is symmetric, we can take $B=\beta T_{k+1}$ for some small enough $\beta>0$ so that the spectral radius of $B$ is strictly less than $1$, satisfying (1), (2), and (4). We may, as needed, also scale $\vec{u}$ to obtain a totally nonzero unit vector $\vec{v}$ in the null space of $B$, satisfying (3).
\end{proof}

Using the matrix $B$ from Lemma \ref{construct} in a construction from Levene et al. from \cite{MR3891770}, we are able to find a matrix $M_{2k+3}$ with the desired properties. Note that \cite{MR3891770} builds a matrix in $\mathcal{S}(\overline{W_{2k+3}})$; our construction differs slightly by using a symmetric matrix $B$, which makes it easier to prove that $M_{2k+3}$ has the SSP.

To describe $M_{2k+3}$ we follow the vertex ordering of $\overline{W_{2k+3}}$ given in \cite{MR3891770}. We label the vertices of $W_{2k+3}$ with the numbers $1,\ldots,2k+3$:  the vertex with degree $k+1$ has label $1$, the vertices with degree $1$ are labelled $2,\ldots, k+2$, and the vertices of degree $2$ are labelled $k+3,\ldots, 2k+3$. So in $W_{2k+3}$, $1$ is adjacent to each of $k+3,\ldots, 2k+3$, and for each $2\leq i\leq k+2$, vertex $i$ is adjacent to $i+k+1$.

\begin{corollary}\label{SSPConst}
For any $k\geq 3$, there is a matrix $M_{2k+3}\in\mathcal{S}(\overline{W_{2k+3}})$ with $q(M_{2k+3})=2$ and the SSP.
\end{corollary}
\begin{proof}
Let $B$ be the $(k+1)\times(k+1)$ matrix and $\vec{v}$ be the vector specified in Lemma \ref{construct}. Using the construction in Equation (6) of \cite{MR3891770} we let
\[
M_{2k+3}=\widehat{M}(B,\vec{v},\alpha)=\left[\begin{array}{rrr}
0 & \vec{v}^T & \vec{0}^T\\
\vec{v} & \sqrt{I_{k+1}-(\alpha^2B^2+\vec{v}\vec{v}^T)} & \alpha B\\
\vec{0} & \alpha B & -\sqrt{I_{k+1}-\alpha^2B^2}
\end{array}\right]=\left[\begin{array}{rrr}
0 & \vec{v}^T & \vec{0}^T\\
\vec{v} & C_1 & \alpha B\\
\vec{0} & \alpha B & C_2
\end{array}\right],
\]
where $\alpha\in [-1,1]$ is chosen as follows. Given that the spectral radius of $B$ is less than $1$, it follows from Lemmas 2.1 and 2.3 in \cite{MR3891770} that $M_{2k+3}$ is a symmetric orthogonal $(2k+3)\times(2k+3)$ matrix, and for all but finitely many values of $\alpha\in[-1,1]$, $C_1,C_2\in\mathcal{S}(K_{k+1})$. Since $B$ has zero diagonal, provided $\alpha\neq 0$, we see that $M_{2k+3}\in\mathcal{S}(\overline{W_{2k+3}})$. We now show that for any such choice of $\alpha$ the resulting $M_{2k+3}$ also has the SSP. We do this by showing directly that the only certificate matrix is $X=O$.

Suppose $X$ is a $(2k+3)\times (2k+3)$ symmetric matrix so that $X\circ M_{2k+3}=X\circ I_{2k+3}=O$ and $XM_{2k+3}=M_{2k+3}X$. From the entrywise product and symmetry constraints, $X$ has the form
\[
X = \left[\begin{array}{rrr}
0 & \vec{0}^T & \vec{u}^T\\
\vec{0} & O & D\\
\vec{u} & D & O
\end{array}\right],
\]
where $\vec{u}$ is unconstrained, and $D$ is a $(k+1)\times (k+1)$ diagonal matrix. Setting
\[
\left[\begin{array}{rrr}
0 & \alpha\vec{u}^TB & \vec{u}^TC_2\\
\vec{0} & \alpha DB & DC_2\\
D\vec{v} & \vec{u}\vec{v}^T+DC_1 & \alpha DB
\end{array}\right]
= XM_{2k+3} = M_{2k+3}X
=\left[\begin{array}{rrr}
0 & \vec{0}^T & \vec{v}^TD\\
\alpha B\vec{u} & \alpha BD & \vec{v}\vec{u}^T+C_1D\\
C_2\vec{u} & C_2D & \alpha BD
\end{array}\right],
\]
we see that equality in the $(2,2)$-entry yields $\alpha DB=\alpha BD$. Since $B$ is totally nonzero off the diagonal, this implies $D=\beta I_{k+1}$ for some number $\beta$.  From the $(2,3)$-entry of $XM_{2k+3}=M_{2k+3}X$ we obtain the equation
\[
DC_2 = \vec{v}\vec{u}^T+C_1D,
\]
which then implies
\[
\beta (C_2-C_1) = \vec{v}\vec{u}^T.
\]

Now we revisit the construction of $M_{2k+3}$ in \cite{MR3891770}. From Lemma 2.3 in  \cite{MR3891770}, we have $C_1=-C_2-\vec{v}\vec{v}^T$. Thus our expression for $\vec{u}\vec{v}^T$ simplifies to
\[
\vec{u}\vec{v}^T= \vec{v}\vec{u}^T = 2\beta C_2+\beta \vec{v}\vec{v}^T.
\]
From the $(2,1)$-entry of $M_{2k+3}X=XM_{2k+3}$ we have $\alpha B\vec{u}=\vec{0}$. Therefore
\[
O = \alpha B\vec{u}\vec{v}^T=2\beta \alpha BC_2+\beta \alpha B\vec{v}\vec{v}^T=2\beta \alpha BC_2
\]
as $\vec{v}$ was chosen so that $B\vec{v}=\vec{0}$.

Note that if $\vec{x}$ is a $\gamma$-eigenvector for $B$, then $\vec{x}$ is a $\sqrt{1-\alpha^2\gamma^2}$-eigenvector for $C_2$ (since the spectral radius of $B$ is strictly less than 1 and $\alpha\in [-1,1]$, $1-\alpha^2\gamma^2> 0$). If $\gamma\neq 0$, then
\[
\vec{0}=O\vec{x}=2\beta \alpha BC_2\vec{x}=-2\beta\alpha\gamma\sqrt{1-\alpha^2\gamma^2}\vec{x},
\]
and $\beta=0$. Since $D=\beta I_{k+1}$, this means that $D=O$. Also since $\vec{u}\vec{v}^T=\beta(C_2-C_1)$ and $\vec{v}$ is totally nonzero, we have that $\vec{u}=\vec{0}$.

Therefore $X=O$ is the only $(2k+3)\times (2k+3)$ symmetric matrix so that $X\circ M_{2k+3}=X\circ I_{2k+3}=O$ and $XM_{2k+3}=M_{2k+3}X$. Thus $M_{2k+3}$ has the SSP, as desired.
\end{proof}

To complete the proof of Theorem \ref{tightbipartitechar}, it remains to consider the original (unsimplified) graph $G$.  By Observation \ref{whysimplified}, if $G$ can be obtained from a simplified graph with $q=2$ by joined-duplication, then $q(G) = 2$.  Any graph $G$ obtained from $P_3$ or $\overline{P_4\cup K_1}$ via a sequence of joined-duplications that maintain $e(\Gbar) = n-2$ has $\overline{G}=S_{a,b}\cup K_1$ with $a\geq 1$ and $b=0$ or $b\geq 1$, respectively. It follows by Lemma \ref{jdup} that $q(G)\leq 3$ for all such graphs.  In the first case, since $G$ contains a unique shortest path $uzv$ where $u$ is a leaf vertex of $S_{a,0}$, $v$ is the non-leaf vertex, and $z$ is the isolated vertex in $\Gbar$, we have $q(G)=3$ by Lemma \ref{nouniquenb}. Similarly, in the second case, $G$ contains a unique shortest path $uzv$ where $u$ and $v$ are the non-leaf vertices of $S_{a,b}$ in $\overline{G}$ and $z$ is the isolated vertex in $\overline{G}$. Therefore we see that $q(G)=3$ for all such graphs.

\section{Further Evidence in Favor of Conjecture \ref{mainconj}}\label{OtherComplements}

From the analysis in the previous sections and appealing to \cite{FallatMojallal}, we have verified Conjecture \ref{mainconj} for $n \leq 8$. Consequently, it is tempting to consider possible inductive strategies that may lead to additional positive evidence for this conjecture. Along these lines we have the following result for general collections of edges removed from the complete graph on $n$ vertices.

\begin{theorem} 
    Suppose $n\geq 4$ and that $H$ is a graph on $n-2$ vertices with $e(H) \leq n-4$. If Conjecture \ref{mainconj} holds for graphs on at most $n-1$ vertices, then $q(K_n-E(H))=2.$
\end{theorem}

\begin{proof}
    Assume, by induction, that Conjecture \ref{mainconj} holds for graphs on at most $n-1$ vertices. Thus, $K_n-E(H)=K_2\vee \overline{H}=K_1\vee (K_1\vee \overline{H})$. Observe that \[e(K_1\vee \overline{H}) \geq\binom{n-1}{2} -(n-4).\]  Hence $q(K_1\vee \overline{H})=2$ by the inductive hypothesis. Then $K_n-E(H)$ can be realized as the joined-duplication of a graph $K_1\vee \overline{H}$ with $q(K_1\vee \overline{H})=2.$ So $q(K_n-E(H))=2.$ 
\end{proof}

Looking beyond bipartite complements, we consider a natural initial case where the subgraph $H$ removed from $K_n$ contains one large odd cycle (relative to $n$), namely $H$ is an odd cycle with $n-3$ edges.   To show $q(\overline{C_{n-3}} \vee K_3)=2$ when $n-3$ is odd, we rely on the construction from Proposition 5.9 of \cite{MR2801365}.  For completeness and to be explicit, we include its proof in the following argument. 

\begin{lemma}\label{MikesLemma}
For $n\geq 6$ there exists an $(n-3)\times 3$ matrix $M$ such that $MM^T \in \mathcal{S}(\overline{C_{n-3}})$ and $M^TM$ has 3 distinct nonzero eigenvalues.
\end{lemma}

\begin{proof}
If $n=6$, the matrix \[M = \begin{bmatrix} 1 & 0 & 0 \\ 0 & 2 & 0\\ 0 & 0 & 3 \end{bmatrix}\]
satisfies the two requirements, so we may assume $n\geq 7$.

Let $\vec{v}_1=[1,1,2]$, $\vec{v}_2=[1,3,-2]$, and $\vec{v}_3=[-1,1,1]$ be row vectors in ${\mathbb R}^3$.  Suppose that nonzero row vectors $\vec v_4, \ldots, \vec v_k \in {\mathbb R}^3$ have been constructed so that 
      \begin{enumerate}
        \item $\vec{v}_i \cdot \vec{v}_{i+1} = 0$ for all $1\leq i \leq k-1$, and
        \item $\vec{v}_i \cdot \vec{v}_j \not=0$ for all $1\leq i,j \leq k$ with $|i-j| > 1$.
    \end{enumerate}
Note that conditions (1) and (2) imply no two of the vectors are scalar multiples of one another.  Consider the set
\[
    S = \mathrm{span}\{\vec{v}_k\}^\perp \setminus \left[ \bigcup_{i=1}^{k-1}  \mathrm{span}\{\vec{v}_i \times \vec{v}_k \}
    \cup
    \bigcup_{i=2}^{k}  \mathrm{span}\{\vec{v}_1 \times \vec{v}_i \}^{\perp}
    \right].
\]
Since $\mathrm{span}\{\vec{v}_i \times \vec{v}_k \}$ is a one-dimensional subspace of $\mathrm{span}\{\vec{v}_k\}^\perp$ for $1\leq i \leq k-1$, and since $\mathrm{span}\{\vec{v}_1 \times \vec{v}_i \}^{\perp}$ must intersect $\mathrm{span}\{\vec{v}_k\}^\perp$ in a one-dimensional subspace for $2\leq i\leq k$ (note $k>2$), we can always choose a nonzero vector $\vec w \in S$.  It follows that $\vec w \cdot \vec v_k = 0$.  Moreover, for $1\leq i \leq k-1$, $\vec w \cdot \vec v_i \neq 0$ since $\vec{w}\notin \mathrm{span}\{\vec{v}_i \times \vec{v}_k \}$.

Let $\vec u = \vec v_1 \times \vec w$.  Then $\vec u \cdot \vec v_i = 0$ if and only if 
\[
0 = (\vec v_1 \times \vec w) \cdot \vec v_i = (\vec v_i \times \vec v_1) \cdot \vec w.
\]
So, for $2\leq i \leq k$, $\vec u \cdot \vec v_i \neq 0$ since $w\notin \mathrm{span}\{\vec{v}_1 \times \vec{v}_i \}^{\perp}$.  Arguing inductively, we may assume $k=n-5$ and $\vec w = \vec v_{k+1}$, so that the vectors $\vec v_1, \ldots, \vec v_{k+1}, \vec u$ are a faithful orthogonal representation of $\overline{C_{n-3}}$.

To build the matrix $M$, let $P$ be the $3\times 3$ matrix whose rows are $\vec v_1, \vec v_2, \vec v_3$, and let $Q$ be the $(n-6)\times 3$ matrix with rows $\vec v_4, \ldots, \vec v_{k+1}, \vec u$.  The eigenvalues of $P^TP$ are $\lambda_1 = 14$, $\lambda_2=7$, and $\lambda_3=2$.  Let $\rho(Q^TQ)$ denote the spectral radius of $Q^TQ$, and find $\varepsilon > 0$ such that $\varepsilon^2 \rho(Q^TQ) < \frac{1}{2}\mathrm{min}\{\lambda_1-\lambda_2, \lambda_2-\lambda_3\}$.  Finally, let
\[ M = \left[ 
\begin{array}{c} P \\  \varepsilon Q \end{array} \right].\]
Then $MM^T \in \mathcal{S}(\overline{C_{n-3}})$ and $M^TM = P^TP + \varepsilon^2 Q^TQ$.  By Weyl's Theorem (see \cite[Corollary 4.3.15]{HJ-new}) 
the spectrum of a symmetric matrix is stable under perturbation; that is, for $j=1,2,3$,
\[ |\lambda_j(M^TM) - \lambda_j(P^TP)| \leq \varepsilon^2 \rho(Q^TQ).\]
Therefore, $M^TM$ has three distinct eigenvalues, as desired.

\end{proof}

\begin{theorem} \label{tri-cyc}
   Suppose $n\geq 6$ is even. Then 
   $q(\overline{C_{n-3}} \vee K_3)=2$, and there is a matrix realization with two distinct eigenvalues and the SSP.
\end{theorem}

\begin{proof}
Let $M$ be the $(n-3)\times 3$ matrix constructed in Lemma \ref{MikesLemma} and let $A = MM^T$. Then $A$ is an $(n-3) \times (n-3)$ positive semidefinite matrix in $\mathcal{S}(\overline{C_{n-3}})$ with rank 3. Furthermore, from the construction of $M$, we have that $q(A)=4$ and that $M^{T}M$ is not a scalar matrix.  Choose $\alpha > 0$ so that $\alpha I_3 - M^T M$ is positive definite and choose a $3 \times 3$ matrix $M_1$ so that $M_1^T M_1=\alpha I_3 - M^T M$.

Let
\[ B = \left[ 
\begin{array}{c} M \\  M_1 \end{array} \right].\]
It follows that $B^TB = \alpha I_3$.

Since $MM_1^T$ may contain zero entries, $BB^T$ may not lie in $\mathcal{S}(\overline{C_{n-3}}\vee K_3)$. So we act on $M_1$, via a $3 \times 3$ orthogonal matrix $R$, to obtain a new matrix
\[ B' = \left[ 
\begin{array}{c} M \\  RM_1\end{array} \right],\] to guarantee that 
\[ B'(B')^T = \left[ \begin{array}{cc}
MM^T & MM_1^T R^T\\ RM_1M^T & RM_1M_1^TR^T
\end{array} \right]\] is in $\mathcal{S}(\overline{C_{n-3}} \vee H)$, where $H$ is some graph on three vertices.  Note that
\[
(B')^TB' = M^T M+ (RM_1)^T RM_1 = \alpha I_3.
\]
By construction, $C=RM_1M_1^T R^T$ has three distinct positive eigenvalues. Thus the graph of $C$ is either empty, contains one edge, two edges, or three edges. It is not difficult to check that if the graph of $C$ contains zero, two, or three edges, then $C$ has the SSP. If the graph of $C$ has only one edge, then it is still true that $C$ will have the SSP. To see this,  suppose the graph of $C$ has only one edge so that $C$ has the form
\[ C= \left[ \begin{array}{ccc}
a & b & 0 \\ b & c & 0 \\ 0 & 0 & d
\end{array} \right], \] and let the certificate matrix be $X = [x_{ij}]$. It follows that $CX$ is symmetric if and only if
\[ \left[ \begin{array}{cc}
a & b  \\ b & c  
\end{array}\right] \left[\begin{array}{c} x_{13} \\ x_{23}\end{array}\right]= d \left[\begin{array}{c} x_{13} \\ x_{23}\end{array}\right].\] 
Hence $d$ is a multiple eigenvalue of $C$, unless $x_{13}=x_{23}=0$, which is a contradiction, as $C$ has three distinct eigenvalues.

Finally, using Lemma 5.1 from \cite{JOS19} the matrix $B'(B')^T$ has the SSP and $q(B'(B')^T)=2$. Thus $q(\overline{C_{n-3}} \vee K_3)=2$ with an SSP realization.
\end{proof}

The proof above can be easily extended to a slightly more general family of graphs, the graphs that can be expressed as $G=\Gamma\vee H$ where $H$ is any graph on three vertices. If $A$ is an $n \times n$ matrix with distinct eigenvalues $\lambda_1$ and $\lambda_2$ with respective multiplicities $m_1$ and $m_2$, then we call $[m_1, m_2]$ the multiplicity bipartition (MB) of $A$.

\begin{porism}
    Let $G =\Gamma\vee \overline{K_3}$, where $\Gamma$ is a graph on $n-3$ vertices with the following properties: $q(\Gamma)=2$ with an SSP realization and with MB $[n-6,3]$. Then $q(G)=2$ with an SSP realization. 
\end{porism}

\begin{proof} 
    Applying the hypotheses on the graph $\Gamma$, we know there exists $B \in \mathcal{S}(\Gamma)$ such that $B=CC^T$, where $C$ is an $(n-3)\times 3$ matrix and $C^T C$ is a nonzero scalar matrix.
    Construct an $n \times 3$ matrix $M$ obtained from $C$ as follows:
\[ M = \left[ 
\begin{array}{c} C \\ Q \end{array} \right],\] 
where $Q$ is orthogonal. Following similar reasoning as in the proof of Theorem \ref{tri-cyc}, using isometries as needed,  we may assume that $MM^T\in\mathcal{S}(G)$. To verify that
\[ MM^T = \left[ \begin{array}{cc}
B & CQ^T \\ QC^T & I_3
\end{array} \right]\]
has the SSP, suppose the matrix $B$ above also has the SSP. Suppose $X$ is a designated real symmetric certificate matrix for the SSP. Then
\[ (MM^T)X =  \left[ \begin{array}{cc}
B & CQ^T \\ QC^T & I_3
\end{array} \right] \left[ \begin{array}{cc}
X_1 & O \\ O & X_2
\end{array} \right]  = \left[ \begin{array}{cc}
BX_1 & (CQ^T)X_2 \\ (QC^T)X_1 & X_2
\end{array} \right].\]
Since $B$ has the SSP, it follows that $X_1=O$. Since $C$ has rank 3, it has no nonzero null vectors, and so $(CQ^T)X_2=O$ implies $Q^T X_2=O$. But $Q$ is orthogonal, so we have $X_2=O$. Hence $MM^T$ has the SSP.

Note that $M^T M = C^TC+Q^T Q$ is a nonzero scalar matrix. Hence $q(G)=2$ with an SSP realization.
\end{proof}

An example class of such graphs $\Gamma$ are complements of trees with the property that: $q(\Gamma)=2$ with  a positive semidefinite SSP realization that has rank equal to 3. A specific example that satisfies the above requirements is 
$\Gamma = \overline{P_n}$, where $n > 4$ and not divisible by 3 (consult the works \cite{MR3891770} for the claim regarding $q$; \cite{AIM} for the claim about the minimum positive semidefinite rank; and 
\cite{JOS19} for guaranteeing an SSP realization).

Theorems \ref{LB1} and \ref{n-3 edges} and the results of this section provide evidence in the affirmative for Conjecture \ref{mainconj}. Theorem \ref{tightbipartitechar} shows that for bipartite complements, the only way to remove $n-2$ edges from a complete graph and obtain $q>2$ is if the edges removed create a unique shortest path of length $2$ between some pair of distinct vertices. In fact, this is the only way we know of to obtain $q>2$ by removing $n-2$ edges even without restricting to removing a bipartite graph. Thus we state the following strengthening of Conjecture \ref{mainconj}.

\begin{conjecture}
    Let $G$ be a graph of order $n\geq 3$.
If $e(\overline{G})\leq  n-2$, then $q(G)=3$ if and only if $\overline{G}= S_{a,b}\cup K_1$ for some $a,b\geq 0$, and $q(G)=2$ otherwise. 
\end{conjecture}

\section*{Acknowledgements}

The authors' collaboration began as part of the ``Inverse Eigenvalue Problems for Graphs and Zero Forcing” Research Community sponsored by the American Institute of Mathematics (AIM). We thank AIM for their support, and we thank the organizers and participants for contributing to this stimulating research experience.  Shaun M. Fallat was supported in part by an NSERC Discovery Research Grant, Application No.: RGPIN--2019--03934.  Veronika Furst and the AIM SQuaRE meetings during which this project was started and completed were supported in part by NSF grant DMS-2331072.  We thank AIM for providing a supportive and mathematically rich environment for our SQuaRE.  Michael Tait was supported in part by NSF grant DMS-2245556 and a Villanova University Summer Grant.

\bibliographystyle{plain}
\bibliography{bibliography}

\begin{thebibliography}{10}

\bibitem{MR4044603}
Mohammad Adm, Shaun Fallat, Karen Meagher, Shahla Nasserasr, Sarah Plosker, and Boting Yang.
\newblock Achievable multiplicity partitions in the inverse eigenvalue problem of a graph.
\newblock {\em Spec. Matrices}, 7:276--290, 2019.

\bibitem{MR3118943}
Bahman Ahmadi, Fatemeh Alinaghipour, Michael~S. Cavers, Shaun Fallat, Karen Meagher, and Shahla Nasserasr.
\newblock Minimum number of distinct eigenvalues of graphs.
\newblock {\em Electron. J. Linear Algebra}, 26:673--691, 2013.
\newblock Erratum available at: \url{https://journals.uwyo.edu/index.php/ela/article/view/1293/5765}.

\bibitem{MR4284782}
John Ahn, Christine Alar, Beth Bjorkman, Steve Butler, Joshua Carlson, Audrey Goodnight, Haley Knox, Casandra Monroe, and Michael~C. Wigal.
\newblock Ordered multiplicity inverse eigenvalue problem for graphs on six vertices.
\newblock {\em Electron. J. Linear Algebra}, 37:316--358, 2021.

\bibitem{MR4074182}
Wayne Barrett, Steve Butler, Shaun~M. Fallat, H.~Tracy Hall, Leslie Hogben, Jephian C.-H. Lin, Bryan~L. Shader, and Michael Young.
\newblock The inverse eigenvalue problem of a graph: multiplicities and minors.
\newblock {\em J. Combin. Theory Ser. B}, 142:276--306, 2020.

\bibitem{MR4603836}
Wayne Barrett, Shaun Fallat, Veronika Furst, Franklin Kenter, Shahla Nasserasr, Brendan Rooney, Michael Tait, and Hein van~der Holst.
\newblock Sparsity of graphs that allow two distinct eigenvalues.
\newblock {\em Linear Algebra Appl.}, 674:377--395, 2023.

\bibitem{MR4646344}
Wayne Barrett, Shaun Fallat, Veronika Furst, Shahla Nasserasr, Brendan Rooney, and Michael Tait.
\newblock Regular graphs of degree at most four that allow two distinct eigenvalues.
\newblock {\em Linear Algebra Appl.}, 679:127--164, 2023.

\bibitem{MR3665573}
Wayne Barrett, Shaun Fallat, H.~Tracy Hall, Leslie Hogben, Jephian C.-H. Lin, and Bryan~L. Shader.
\newblock Generalizations of the strong {A}rnold property and the minimum number of distinct eigenvalues of a graph.
\newblock {\em Electron. J. Combin.}, 24(2):Paper No. 2.40, 28, 2017.

\bibitem{MR3904092}
Beth Bjorkman, Leslie Hogben, Scarlitte Ponce, Carolyn Reinhart, and Theodore Tranel.
\newblock Applications of analysis to the determination of the minimum number of distinct eigenvalues of a graph.
\newblock {\em Pure Appl. Funct. Anal.}, 3(4):537--563, 2018.

\bibitem{MR2801365}
Matthew Booth, Philip Hackney, Benjamin Harris, Charles~R. Johnson, Margaret Lay, Terry~D. Lenker, Lon~H. Mitchell, Sivaram~K. Narayan, Amanda Pascoe, and Brian~D. Sutton.
\newblock On the minimum semidefinite rank of a simple graph.
\newblock {\em Linear Multilinear Algebra}, 59(5):483--506, 2011.

\bibitem{FallatMojallal}
Shaun Fallat and Seyed~Ahmad Mojallal.
\newblock Spectral applications of vertex-clique incidence matrices associated with a graph.
\newblock {\em Mathematics}, 11(16), 2023.

\bibitem{MR3240026}
Cheryl Grood, Johannes Harmse, Leslie Hogben, Thomas~J. Hunter, Bonnie Jacob, Andrew Klimas, and Sharon McCathern.
\newblock Minimum rank with zero diagonal.
\newblock {\em Electron. J. Linear Algebra}, 27:458--477, 2014.

\bibitem{AIM}
AIM Minimum Rank-Special Graphs~Work Group et~al.
\newblock Zero forcing sets and the minimum rank of graphs.
\newblock {\em Linear Algebra Appl.}, 428(7):1628--1648, 2008.

\bibitem{MR3506498}
Keivan Hassani~Monfared and Bryan~L. Shader.
\newblock The nowhere-zero eigenbasis problem for a graph.
\newblock {\em Linear Algebra Appl.}, 505:296--312, 2016.

\bibitem{MR4478249}
Leslie Hogben, Jephian C.-H. Lin, and Bryan~L. Shader.
\newblock {\em Inverse problems and zero forcing for graphs}, volume 270 of {\em Mathematical Surveys and Monographs}.
\newblock American Mathematical Society, Providence, RI, 2022.

\bibitem{HJ-new}
Roger~A. Horn and Charles~R. Johnson.
\newblock {\em Matrix analysis}.
\newblock Cambridge University Press, Cambridge, second edition, 2013.

\bibitem{MR3891770}
Rupert~H. Levene, Polona Oblak, and Helena \v{S}migoc.
\newblock A {N}ordhaus-{G}addum conjecture for the minimum number of distinct eigenvalues of a graph.
\newblock {\em Linear Algebra Appl.}, 564:236--263, 2019.

\bibitem{MR4774297}
Rupert~H. Levene, Polona Oblak, and Helena \v{S}migoc.
\newblock Distinct eigenvalues are realizable with generic eigenvectors.
\newblock {\em Linear Multilinear Algebra}, 72(12):2054--2068, 2024.

\bibitem{JOS19}
Jephian C.-H Lin, Polona Oblak, and Helena Šmigoc.
\newblock The strong spectral property for graphs.
\newblock {\em arXiv preprint}, 2019.

\end{thebibliography}

\end{document}